\documentclass{amsart}

\address[hirachi@ms.u-tokyo.ac.jp]{
Graduate School of Mathematical Sciences, 
The University of Tokyo,
3-8-1 Komaba, Meguro, Tokyo 153-8914 JAPAN}

\newtheorem{theorem}{Theorem}[section]

\theoremstyle{definition}

\newcommand{\calC}{\mathcal{C}}

\newcommand{\calE}{\mathcal{E}}

\newcommand{\calH}{\mathcal{H}}

\newcommand{\calJ}{\mathcal{J}}

\newcommand{\calM}{\mathcal{M}}
\newcommand{\calN}{\mathcal{N}}
\newcommand{\calO}{\mathcal{O}}

\newcommand{\calQ}{\mathcal{Q}}
\newcommand{\calR}{\mathcal{R}}

\newcommand{\bC}{\mathbb{C}}

\newcommand{\bP}{\mathbb{P}}
\newcommand{\bR}{\mathbb{R}}

\renewcommand{\Re}{\operatorname{Re}}
\renewcommand{\Im}{\operatorname{Im}}
\newcommand{\Ric}{\operatorname{Ric}}

\newcommand{\Scal}{\operatorname{Scal}}

\newcommand{\vol}{\operatorname{vol}}

\newcommand{\wt}{\widetilde}
\newcommand{\wf}{\wt f}

\newcommand{\pa}{\partial}
\newcommand{\wh}{\widehat}
\renewcommand{\th}{\theta}
\newcommand{\up}{\Upsilon}
\newcommand{\conj}{\overline}

\newcommand{\comment}[1]{\relax}

\newcommand{\abs}[1]{\lvert#1\rvert}

\renewcommand\th{\theta}

\newcommand\ol{\overline}

\newcommand\YDsym{\setlength{\unitlength}{0.25mm}
\begin{picture}(21,11)(0,1)
\put(0,0){\line(0,1){10}}
\put(10,0){\line(0,1){10}}
\put(20,0){\line(0,1){10}}
\put(0,10){\line(1,0){20}}
\put(0,0){\line(1,0){20}}
\end{picture}}

\newcommand\YDsymsc{{\setlength{\unitlength}{0.15mm}
\begin{picture}(20,10)(-1,2)
\put(0,0){\line(0,1){10}}
\put(10,0){\line(0,1){10}}
\put(20,0){\line(0,1){10}}
\put(0,10){\line(1,0){20}}
\put(0,0){\line(1,0){20}}
\end{picture}}}

\newcommand\YDN{\setlength{\unitlength}{0.25mm}
\begin{picture}(20,21)(0,6)
\put(0,0){\line(0,1){20}}
\put(10,0){\line(0,1){20}}
\put(20,10){\line(0,1){10}}
\put(0,10){\line(1,0){20}}
\put(0,20){\line(1,0){20}}
\put(0,0){\line(1,0){10}}
\end{picture}}

\newcommand\YDND{\setlength{\unitlength}{0.25mm}
\begin{picture}(20,21)(0,0)
\put(0,0){\line(0,1){20}}
\put(10,0){\line(0,1){20}}
\put(20,10){\line(0,1){10}}
\put(0,10){\line(1,0){20}}
\put(0,20){\line(1,0){20}}
\put(0,0){\line(1,0){10}}
\end{picture}}

\newcommand\YDbox{
\setlength{\unitlength}{0.25mm}
\begin{picture}(21,20)(0,0)
\put(0,0){\line(0,1){20}}
\put(10,0){\line(0,1){20}}
\put(20,0){\line(0,1){20}}
\put(0,0){\line(1,0){20}}
\put(0,10){\line(1,0){20}}
\put(0,20){\line(1,0){20}}
\end{picture}}

\newcommand\YDboxsmall{
\setlength{\unitlength}{0.17mm}
\begin{picture}(21,20)(0,3)
\put(0,0){\line(0,1){20}}
\put(10,0){\line(0,1){20}}
\put(20,0){\line(0,1){20}}
\put(0,0){\line(1,0){20}}
\put(0,10){\line(1,0){20}}
\put(0,20){\line(1,0){20}}
\end{picture}}

\newcommand\YDs{
\setlength{\unitlength}{0.25mm}
\begin{picture}(10,10)(0,0)
\put(0,0){\line(0,1){10}}
\put(10,0){\line(0,1){10}}
\put(0,0){\line(1,0){10}}
\put(0,10){\line(1,0){10}}
\end{picture}}

\newcommand\YDmid{
\setlength{\unitlength}{0.25mm}
\begin{picture}(20,20)(0,0)
\put(0,-10){\line(0,1){30}}
\put(10,-10){\line(0,1){30}}
\put(20,0){\line(0,1){20}}
\put(0,0){\line(1,0){20}}
\put(0,-10){\line(1,0){10}}
\put(0,10){\line(1,0){20}}
\put(0,20){\line(1,0){20}}
\end{picture}}

\title[$Q$ and $Q$-prime curvature]{%
$Q$ and $Q$-prime curvature in CR geometry
}

\author[Kengo Hirachi]
{Kengo Hirachi
}

\begin{document}

\begin{abstract}
The $Q$-curvature has been playing a central role in conformal geometry since its discovery by T. Branson.
It has natural analogy in CR geometry, however, the CR $Q$-curvature vanishes on the boundary of a strictly pseudoconvex domain in $\bC^{n+1}$ with a natural choice of contact form.  This fact enables us to define a ``secondary'' $Q$-curvature,
which we call $Q$-prime curvature (it was first introduced by  J. Case and P. Yang in the case $n=1$).  The integral of the $Q$-prime curvature, the total $Q$-prime curvature,  is a CR invariant of the boundary.  When $n=1$, it agrees with the Burns-Epstein invariant, which is a Chern-Simons type invariant in CR geometry.  For all $n\ge1$, it has non-trivial variation under the deformation of domains.  Combining the variational formula with the  deformation complex of CR structures, we show that the total $Q$-prime curvature takes local maximum at the standard CR sphere in a formal sense.

This talk is a report in collaboration with  Rod Gover, Yoshihiko Matsumoto, Taiji Marugame
and Bent \O{}rsted.
\end{abstract}

\maketitle

\section{Introduction}
In 1979, C. Fefferman \cite{F2} proposed a program of studying the analysis and geometry of strictly pseudoconvex domains in $\bC^{n+1}$, which is called ``Parabolic invariant theory''.
The basic idea is to consider the asymptotic expansions of the Bergman and Szeg\"o kernels as analogies of the heat kernel on Riemannian manifolds and develop invariant theory that leads to index theorems.  Since then there has been a series of works that completes the local theory of the asymptotic expansion; see \cite{H0} for an overview.  The main tool is the ambient metric, a Ricci-flat Lorentz-K\"ahler metric defined on the canonical bundle of the domain; all local invariants of the CR structure on the boundary of the domain can be expressed as invariant polynomials in the jets of the curvature of the ambient metric.  However, the connection of the local formula to the global invariants has not been understood.  In this talk, as a continuation of the program, we will give a construction of global CR invariants which we call total $Q$ and  $Q$-prime curvatures.

The concept of $Q$-curvature was first introduced by T. Branson in his study of the functional determinant of conformally invariant differential operators.  In conformal geometry, it is a volume form valued local Riemannian invariant $Q(g)$ defined on even dimensional Riemannian manifolds $(M,g)$.
While $Q$ is not a local conformal invariant, its de Rham class is conformally invariant, i.e.,  the integral $\overline{Q}=\int Q$ over $M$ is a global conformal invariant, which is called the total $Q$-curvature.
The total $Q$-curvature got much attention by the seminal work of Graham-Zworski \cite{GZ} on the scattering theory of conformally compact Einstein manifolds $(X,g_+)$.  The $Q$-curvature on the boundary at  infinity $M=\pa X$ can be expressed in terms of scattering-matrix, and $\overline{Q}$ is characterized as a coefficient in the volume renormalization of $(X,g_+)$.  This approach has natural connection with AdS/CFT correspondence in string theory 
and  gives new insights in conformal geometry.  A nice overview of these progress is given in the lecture notes by Juhl \cite{BJ}.

The $Q$-curvature in CR geometry was first introduced in \cite{H1} in  the description of the Szeg\"o kernel for domains in $\bC^{2}$; the general definition was given later in \cite{FH}. 
See also the results on the Szeg\"o kernel along this line \cite{H3, Bo}.  However, it also turned out that the total $Q$-curvature always vanishes for the boundaries of strictly pseudoconvex domains in $\bC^{n+1}$.   A breakthrough was made recently by Case-Yang \cite{CY}; they defined a ``secondary'' $Q$-curvature $Q'$ on  3-dimensional CR manifolds and called it  $Q$-prime curvature. It is shown that the integral $\overline Q'=\int Q'$, the total $Q$-prime curvature, agrees with the Burns-Epstein invariant \cite{BE1}, which is a secondary invariant in 3-dimensional CR geometry analogous to the Chern-Simons invariant of conformal $4$-manifolds.
 The generalization of $Q'$ to higher dimensions was given by \cite{H4}. When $n\ge2$, $\overline Q'$ is different from the secondary characteristic number defined by Burns-Epstein \cite{BE2}.
The crucial property of $\overline Q'$ is its variational formula under the deformation of domains, which enables us to study  $\overline Q'$ as a functional on the moduli space of CR structures.
In particular, we show that the critical points are exactly the boundaries of the domains with smooth solution to a Dirichlet problem for the complex Monge-Amp\`ere equation  (it can be also characterized as the vanishing of obstruction function, which is a local CR invariant of the boundary).

We can find an intimate relation between the variational formula of $\overline Q'$ and the deformation complex of CR structures as an application of representation theory.
In the parabolic invariant theory, CR manifolds are seen as geometric structures modeled on the homogeneous space $G/P$, where $G=SU(n+1,1)$ and $P$ is a parabolic subgroup.  The invariant theory of $P$ gives much information of the geometry.  The deformation complex of CR structures is the resolution of the adjoint representation $\frak{g}$ by the complex of CR invariant (or equivalently $G$-invariant) differential operators acting on the bundles induced from irreducible representations of $P$; this is known as an example of (generalized)  Bernstein-Gelfand-Gelfand complex \cite{CSS}. 
 The Hessian of $\overline Q'$ at the sphere $S^{2n+1}=G/P$ gives a CR invariant, self-adjoint, differential operator $L$ between two bundles in the deformation complex and the CR invariance property of $L$ forces the operator to be semidefinite.
 Moreover, the kernel of $L$ can be described by using the operators in the complex. Thus we can read  geometric information from the Hessian.  As a result we conclude that $\overline{Q}'$ takes local maximum at the standard CR sphere in a formal sense; see Theorem \ref{mainthm}.
  We  also apply the same procedure to the total $Q$-curvature for partially integrable CR manifolds (which may not be embeddable).  In this generalized setting,
the total $Q$-curvature is non-trivial and the Hessian has an interesting connection to the integrability condition.

The computation of the Hessian of $\overline Q'$ is an analogy of that of $\overline Q$ in conformal geometry given in M\o{}ller-\O{}rsted
\cite{MO}.  In the conformal case, the deformation complex is simpler and it is easy to describe the kernel of the Hessian.  We include an overview of the conformal $Q$ in the next section, which should help to understand the more involved structure of CR $Q$ and $Q'$.

\section{$Q$-curvature in conformal geometry}
We start with a quick review of $Q$-curvature in the conformal geometry with an intention to explain the relation between the total $Q$ and the deformation complex.  Many deep results in geometric analysis of $Q$-curvature are not mentioned; the article by Alice Chang \cite{AC} gives a clear overview in this direction.

 \
\subsection{Dimensions 2 and 4}
The $Q$-curvature is defined as a generalization of the scalar curvature on surfaces to higher even dimensions in the context of conformal geometry.
Given a Riemannian manifold $(M,g)$ of dimension $n$,
let $\Scal$ be the scalar curvature and $\Delta=-g^{ab}\nabla_{a}\nabla_{b}=-\nabla_{a}\nabla^{a}$ be the Laplacian (we use the Einstein summation convention). 
 If we denote by $\wh\Scal$ and $\wh\Delta$ the ones for the scaled metric $\wh g=e^{2\Upsilon}g$ for $\Upsilon\in  C^\infty(M)$, then, for $n=2$, we have
$$
\wh\Scal=e^{-2\Upsilon}(\Scal+2\Delta\Upsilon),\quad \wh\Delta=e^{-2\Upsilon}\Delta.
$$
The factor $e^{-2\Upsilon}$ can be considered as the scaling of the volume form $d\vol_g$.
So  setting
$$
Q_2=\frac12\Scal \cdot d\vol_g,\qquad
P_2 f =\Delta f\cdot d\vol_{g},
$$
we may write the transformation rules as
$$
\wh Q_{2}=Q_2+P_{2}\Upsilon,\qquad \wh P_{2}=P_{2}. 
$$
If $M$ is compact, then $\int_M P_{2}\Upsilon =0$ and thus $\int_{M} Q_{2}$ is a conformal invariant. In fact, it is topological 
as  the Gauss-Bonnet theorem gives
 $\int_{M} Q_2=2\pi\chi(M)$.

On a conformal manifold $(M,[g])$ of even dimension $n$,
the $Q$-curvature $Q_{n}$ and conformally invariant differential operator $P_n$
 are defined as the pair of objects that generalize these properties.  Namely,
$Q_{n}$ is a volume form valued local invariant of metric $g$ 
such that
$$
\wh Q_n=Q_n+P_n\Upsilon
$$
with a conformally invariant differential operator
$$
P_n\colon C^\infty(M)\to C^\infty(M,\wedge^nT^*M),
$$
which is self-adjoint and $P_{n}1=0$. 
The last two properties of $P_{n}$ ensure that the integral
$$
\overline Q_n=\int_M Q_n
$$
is a global invariant of a conformal manifold $(M,[g])$,
which is called the {\em total $Q$-curvature}.
 As we see below, 
  when $n\ge4$, $\overline Q_n$ gives an interesting invariant which is not just topological.

In the case $n=4$, we can define $Q_{4}$ and $P_{4}$ by explicit formulas:
$$
\begin{aligned}
Q_4&=\left(\frac16 \Delta\Scal-\frac12 \Ric^{ab}\Ric_{ab}+\frac16\Scal^2
\right)d\vol,
\\
P_4 f&=\left( \Delta^{2}f+\nabla_{a}(2 \Ric^{ab}-\frac23\Scal g^{ab})\nabla_{b}f\right)d\vol,
\end{aligned}
$$
where $\Ric_{ab}$ is the Ricci tensor.
We may also write $Q_{4}$ as
$$
Q_4=2\operatorname{Pfaff_4}-\frac14|\operatorname{Weyl}|^2d\vol+
\frac16\Delta\Scal\cdot d\vol,
$$
where $\operatorname{Pfaff}_n$ is the Pfaffian, which integrates to $(-2\pi)^{n/2}\chi(M)$, and $|\operatorname{Weyl}|^2$  is the squared norm of the Weyl curvature $\operatorname{Weyl}_{abcd}$, the trace-free part of Riemannian curvature $R_{abcd}$.  Hence the total $Q$-curvature satisfies
$$
\overline{Q}_4=8\pi^2\chi(M)-\frac14\int_M|\operatorname{Weyl}|^2 d\vol.
$$
Since $|\operatorname{Weyl}|^2 d\vol$ is independent of the scale,  the second term is conformally invariant itself.
In particular, we see that 
$$
\overline{Q}_4\le 8\pi^2\chi(M)
$$
and the equality holds if and only if $[g]$ is conformally flat.

\subsection{The ambient metric}  To construct $Q_{n}$ and $P_{n}$ in higher dimensions, we 
 use the ambient metric of Fefferman and Graham \cite{FG2}.  To motivate the definition of the ambient metric, we first recall the M\"obius transformations of the standard sphere and associated metrics.

Let $G=SO(n+1,1)$ be the orthogonal group for the quadratic form
$$
B(\zeta)=-\zeta_0^2+\zeta_1^2+\cdots+\zeta_{n+1}^2,
\qquad
(\zeta_0,\dots,\zeta_{n+1})\in\bR^{n+2}.
$$
Then $G$ preserves the light cone $\calN=\{\zeta\in\bR^{n+2}\setminus\{0\}:B(\zeta)=0\}$ and the hyperboloid $\calH=\{\zeta\in\bR^{n+2}:B(\zeta)=-1\}$.  
The projectivization of $\calN$ can be identified with the unit sphere $S^n=\{x\in\bR^{n+1}:|x|^2=1\}$ by 
$$
S^n\ni x\mapsto \bR(1,x)\in \calN/\bR^*\subset \bP^{n+1}
$$
and $\calH$ can be identified with the unit ball $B^{n+1}=\{x\in\bR^{n+1}:|x|^2<1\}$ by 
$$
 \calC\ni(\zeta_0,\zeta')\mapsto \frac{\zeta'}{1+\zeta_0}\in B^{n+1}.
$$
The action of $G$  on  $\bR^{n+2}$ can be also seen as the isometries of the Lorentzian metric:
$$
\wt g=-d\zeta^2_0+d\zeta_1^2+\cdots+d\zeta_{n+1}^2.
$$
Since $\wt g$ induces the Poincar\'e  metric $g_+$ on $\calH\cong B_{n+1}$,
$$
g_+=4\frac{dx_1^2+\cdots+dx_{n+1}^2}{(1-|x|^2)^2},
$$
we see that the action of $G$ on $B_{n+1}$ gives the isometries of $g_{+}$.

On the other hand, $\wt g$ induces a degenerate two tensor on $\calN$.
For each section of $\pi\colon\calN\to\calN/\bR^*\cong S^n$,
the pullback of $\wt g$ gives a Riemannian metric which is conformal to the standard metric $g_0$ on $S^n$.  Thus (upper half of) $\calN$ can be identified with the metric bundle over $(S^n,[g_0])$; hence $G$ acts as conformal maps of $S^{n}$.

To sum up, we have three spaces of different dimensions on which $G$ acts as automorphisms:
\begin{itemize}
\item Lorentzian space $(\bR^{n+2},\wt g)$;
\item Poincar\'e ball $(B^{n+1}, g_+)$;
\item Conformal sphere $(S^n,[g_0])$.
\end{itemize}

Fefferman-Graham \cite{FG2} generalized these correspondences to a curved conformal manifold $(M,[g])$.  They call $\wt g$ and $g_+$, respectively, the {\em ambient metric} and the {\em Poincar\'e-Einstein metric}.
We first recall the construction of the ambient metric.
Let $\wt M=\bR_+\times M\times (-1,1)$ and choose a local coordinate system $(t,x,\rho)$.  
For each scale $g\in[g]$, the trivial $\bR_+$-bundle $\calM=\bR_+\times M$ can be identified with the metric bundle by $(t,x)\mapsto t^2 g(x)\in S^2T^*M$ and we embed $\calM$ into $\wt M=\calM\times(-1,1)$ as a
hypersurface $\rho=0$.
On $\wt M$, we  consider a Lorentzian metric  of the form
$$
\wt g=2tdt\, d\rho+2\rho dt \,dt
+t^2g_{ab}(x,\rho)dx^a dx^b,
$$
where $g(x,\rho)=g_{ab}(x,\rho)dx^a dx^b$ is a family of Riemannian metrics on $M$ with parameter $\rho$ such that  $g(\,\cdot\,,0)=g$.
To fix $\wt g$, we impose the Einstein equation along $\rho=0$: 
\begin{equation}\label{Ricci-flat-eq}
\Ric(\wt g)=\begin{cases}
O(\rho^\infty) & n$\text{ odd},$\\
O^+(\rho^{n/2-1})
& n$\text{ even.}$
\end{cases}
\end{equation}
Here $f=O^+(\rho^l)$ means that each component of $f$ is $O(\rho^l)$ and $\rho^{-l} f|_{T\calM}$ in the frame 
$dt, dx^a$ is of the form
$$
\begin{pmatrix}
0 & 0 \\
0&\phi_{ab}
\end{pmatrix}
\quad\text{with } g^{ab} \phi_{ab}=0.
$$
The ambient metric $\wt g$ is defined as the solution to the equation \eqref{Ricci-flat-eq}, which exists uniquely modulo $O(\rho^\infty)$ for odd $n$
and modulo $O^+(\rho^{n/2})$ for even $n$.
For even $n$, we set
$$
c_n\big((-\rho)^{1-n/2}\Ric(\wt g)\big)|_{T\calM}=\begin{pmatrix}
0 & 0 \\
0&\calO_{ab}
\end{pmatrix}, \quad
c_n={2^{n-3}(n/2-1)!^2}.
$$
The tensor $\calO_{ab}$ is called the {\em obstruction tensor} and is shown to be a local conformal invariant
in the sense that $\wh\calO_{ab}=e^{(2-n)\Upsilon}\calO_{ab}$
under the change of scale $\wh g=e^{2\Upsilon}g$ that gives the initial data on $\rho=0$.  
We can show that $\Ric(\wt g)=O(\rho^\infty)$ admits a smooth solution if and only 
 if $\calO_{ab}=0$, as its name suggests.
 
While this definition of $\wt g$ depends on a choice of scale $g$, the ambient metric is determined by the conformal class $[g]$ in the following sense: if $\wt g$ and $\wt g'$ are ambient metrics determined by $g,g'\in[g]$, then there is a diffeomorphism
$\Phi(t,x,\rho)=(t',x',\rho')$ such that $\Phi(\lambda t,x,\rho)=(\lambda t',x',\rho')$ for $\lambda>0$,
$\Phi(\calM)=\calM$, and $\Phi^{*}\wt g'=\wt g$  modulo $O(\rho^\infty)$ for odd $n$
and modulo $O^+(\rho^{n/2})$ for even $n$.

The  Poincar\'e-Einstein metric $g_+$  on
$X=M\times (0,1)\ni(x,r)$ is then defined by the pullback of $\wt g$ by the embedding
$X\hookrightarrow\wt M$, 
$(x,r)\mapsto(1/r,x,-r^2/2)$. In the coordinate system $(x,r)$, we have
\begin{equation}\label{g-plus}
g_+=\frac{dr^2+h_{ab}(x,r)dx^{a}dx^{b}}{r^{2}},
\end{equation}
where $h_{ab}(x,r)=g_{ab}(x,-r^2/2)$.  We call $r$ the defining function normalized by the scale $h(\,\cdot\,,0)=g\in [g]$.
The equation \eqref{Ricci-flat-eq} then implies
$$
\Ric(g_+)+n g_+=\begin{cases}
O(r^\infty) & n$\text{ odd}$,\\
c_n^{-1}\calO_{ab} r^{n-2}+O(r^{n-1})
& n$\text{ even}$.
\end{cases}
$$
Alternatively, one may define Poincar\'e-Einstein metric to be the solution to this equation of the form \eqref{g-plus}.

\subsection{Total $Q$-curvature and volume renormalization}
Now we are ready to define $Q_n$ and $P_n$.
For a conformal manifold $(M,[g])$ of even-dimension $n$,
we take the ambient metric $\wt g$.  Let $\wt\Delta=-\wt\nabla^A\wt\nabla_A$ be the (wave) Laplacian; here $\wt \nabla$ is the Levi-Civita connection of $\wt g$ and the index $A$
runs through $0,1,\dots, n+1$.
The $Q$-curvature and the invariant operator $P_{n}$  for the scale $g\in[g]$ are given by
$$
\begin{aligned}
Q_n&=-(\wt\Delta^{n/2}\log t)\big|_{t=1,\rho=0}d\vol_g,
\\
P_nf&=(\wt\Delta^{n/2} f)\big|_{t=1,\rho=0}d\vol_g.
\end{aligned}
$$
Here $f\in C^\infty(M)$ is identified with a function on 
$\wt M$ that is free of $(t,\rho)$.  It is shown that
 $P_n$ is independent of the choice of the scale $g\in[g]$ and the pair 
$Q_{n}$, $P_{n}$ satisfies the required properties \cite{FH}.
(See also \S3.2 for the original definition by  T. Branson.)

To study $\overline Q_{n}$, it is useful to give its characterization in
 terms of the Poincar\'e-Einstein metric $g_{+}$.
The complete Riemannian manifold $(X,g_+)$ has infinite volume; we define its finite part by considering the volume expansion of the sub\-domains $\{r>\epsilon\}$:
$$
\begin{aligned}
\int_{r>\epsilon}d\vol_{g_+}=
a_0\epsilon^{-n}+&a_2\epsilon^{-n+2}+\cdots
\\
&+\begin{cases}
a_{n-1}\epsilon^{-1}+V+o(1) & n\text{ odd},\\
a_{n-2}\epsilon^{-2}-L\log\epsilon+V+o(1)
& n\text{ even}.
\end{cases}
\end{aligned}
$$
Here $r$ is a defining function of $M$ normalized by a scale $g\in[g]$ and $a_{j}, V, L$ are constants.  The constant term $V(g_{+},g)$ is called the {\em renormalized volume} of $g_+$ with respect to the scale $g$.
\begin{theorem}[Graham-Zworski \cite{GZ}]
For odd $n$, $V$ is independent of the choice of a scale $g\in[g]$.
For even $n=2m$, $L$ is independent of the choice of a scale $g$.
Moreover, one has
$$
\overline{Q}_n=(-1)^m2^{n-1}m!(m-1)!L.
$$
\end{theorem}

The proof is based on the scattering theory of $\Delta_+$, the Laplacian of $g_+$.  A simpler proof using the Dirichlet problem for $\Delta_{+}$ was later given by Fefferman-Graham \cite{FG1}.

For even $n$, $V(g_+,g)$ is not conformally invariant.
In fact, the variation of the scale gives
$$
\frac{d}{dt}\Big|_{t=0}V(g_+,e^{2t\Upsilon}g)=\int_M\Upsilon v_n(g)d\vol_{g},
$$
where $v_n(g)$ is a local invariant of $g$ called the {\em holographic anomaly}. 
It appears in the expansion of the volume form:
$$
\frac{\vol (h(\cdot,r))}{\vol(g)}=1 +v_2 r^2+v_4 r^4+\cdots.
$$
It is not difficult to see 
$
L=\int_M v_n(g)d\vol_{g}.
$
So $v_n$ has a similar property to $Q_n$, while the transformation law of $v_n$ under the scaling is not easy to write down. In dimension $4$, one has
$$
Q_4=16 v_4+\frac16\Delta\Scal.
$$
General relation between $Q_n$ and $v_j$ has been studied extensively by Juhl and  Fefferman-Graham; see \cite{BJ} and
\cite{FG3}.

We next study the variation of $\overline Q_n$ under the deformation of conformal structures.
Let $g^t$ be a one parameter family of Riemannian metrics on $M$ that preserves the volume form.  Then it has an expansion
$$
g^t_{ab}=g_{ab}+t\psi_{ab}+O(t^2)
$$
with $g^{ab}\psi_{ab}=0$.  We shall denote the total $Q$-curvature for the conformal structure $[g^{t}]$ by $\overline Q_{n}(g^t)$. 

\begin{theorem}[\cite{GH}] The first variation of the total $Q$-curvature is given by
$$
\frac{d}{dt}\Big|_{t=0}{\overline Q}_{n}(g^t)=(-1)^{n/2}
\int_M  \calO_{ab}\psi^{ab}\, d\vol_{g},
$$
where $\calO_{ab}$ is the obstruction tensor for $g=g^{0}$.
\end{theorem}

The critical points of the functional $\overline Q_n$ on the space of conformal structures  are characterized by $\calO_{ab}=0$.
At a critical point, the second variation is given by
$$
\frac{d^2}{dt^2}\Big|_{t=0}{\overline Q}_{n}(g^t)=
\int_M  L_n(\psi)_{ab}\psi^{ab}\, d\vol_{g},
$$
where $L_n$ is a conformally invariant differential operator of order $n$:
$$
L_{n}(\psi)_{ab}=(-1)^{n/2}\frac{d}{dt}\Big|_{t=0}\calO_{ab}(g^t).
$$
By analyzing the eigenvalues of $L_n$, one obtain the following rigidity of the conformal sphere.

\begin{theorem}[M\o ller-\O rsted \cite{MO}]
The total $Q$-curvature has a local maximum at the standard conformal sphere  $(S^n,[g_0])$, $n\ge4$.  Namely, there exists a neighborhood $U$ of $g_0$ in the space of Riemannian metrics on $S^n$ such that
$$
\overline Q_n(g)\le\overline Q_n(g_0)\quad\text{for all }g\in U
$$
and the equality holds if and only if $(S^{n},[g])$ is conformally equivalent to $(S^n,[g_0])$.\end{theorem}

The proof is based on an analysis of  $L_{n}$ using representation theory.  To formulate it, let us recall the deformation complex of conformal structures.
Let $D_0$ be the conformal Killing map
$$D_{0}\colon \Gamma(T^*M)\to\Gamma(S^{2}_{0}T^*M)$$
defined by
$$
D_0( f_a)= \text{trace-free part of }\nabla_{(a}f_{b)},
$$
where
$S^{2}_{0}T^*M$ denote the bundle of symmetric trace-free 2-tensors.
Let
$$D_1\colon \Gamma(S_{0}^{2}T^*M)\to\Gamma(\otimes^{4}T^*M)
$$
 be 
 the linearization of the Weyl curvature:
$$
D_1(\varphi_{ab})=\text{projection to $\YDboxsmall_{0}$ part of }\nabla_{ab}\varphi_{cd}.
$$
Here $\YDboxsmall_{0}$ denotes  the space of trace-free $4$-tensor with symmetry given by the Young diagram $\YDboxsmall$.  Then we have $\Im D_{0}=\ker D_{1}$ and these maps can be extended to a complex, which is known as the deformation complex.  On $S^{6}$, it is given by (thick long arrow omitted)
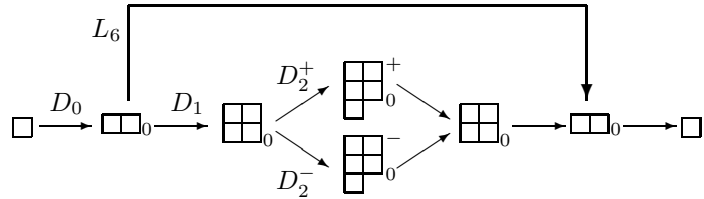
\begin{equation*}\label{deform}
\setlength{\unitlength}{0.7mm}
\begin{picture}(130,37)(-10,-5)
\put(5,24){$L_6$}
{\thicklines
\put(12,12){\line(0,1){18}}
\put(12,30){\line(1,0){87}}
\put(99,30){\vector(0,-1){18}}
}
\put(-10,5){$\YDs$}
\put(-3,10){$D_0$}
\put(-5,7){\vector(1,0){10}}
\put(7,6){$\YDsym_0$}
\put(20,10){$D_1$}
\put(17,7){\vector(1,0){10}}
\put(30,4){$\YDbox_0$}
\put(40,8){\vector(3,2){10}}
\put(40,6){\vector(3,-2){10}}
\put(53,12){$\YDmid_{\, 0}^{\, +}$}
\put(53,-2){$\YDmid_{\, 0}^{\, -}$}
\put(40,-5){$D_2^-$}
\put(40,15){$D_2^+$}
\put(63,15){\vector(3,-2){10}}
\put(63,-1){\vector(3,2){10}}
\put(75,4){$
\YDbox_0$}
\put(96,6){$
\YDsym_0$}
\put(85,7){\vector(1,0){10}}
\put(106,7){\vector(1,0){10}}
\put(117,5){$\YDs$}
\end{picture}
\end{equation*}
\newpage
Here, we only write the symmetries of the tensor bundles in terms of the Young diagram; the superscripts $\pm$ denote the self-dual/anti-self-dual parts.  If we properly put a density weight on each bundle, this becomes a complex of $G$-invariant differential operators. (Recall that $G=SO(n+1,1)$, which acts on $S^{n}$ as conformal maps.)
For general dimensions, the deformation complex has length $n+1$, like the de Rham complex.  The de Rham and deformation complexes are examples of generalized  Bernstein-Gelfand-Gelfand (BGG) complexes in the parabolic geometry modeled on $G/P$.
Each BGG complex gives a resolution of a finite dimensional irreducible representation of $G$; the de Rham complex (resp.~deformation complex)  corresponds to the trivial representation $\bR$ (resp.~the adjoint representation $\frak{g}$). See \cite{CSS}.  

 The Hessian $L_{n}$  gives a $G$-invariant operator from $\YDsym_0$ on the left to $\YDsym_0$ on the right.  Such an operator is unique up to a constant multiple and turns out to be a semi-definite, self-adjoint, operator with kernel
$\Im D_0$.  Since $\Im D_{0}=\ker D_1$ is the tangent space to the submanifold consisting of flat conformal structures on $S^{n}$, we get an infinitesimal version of the theorem.  This analysis can be applied to many other conformal functionals, e.g., the determinant of Yamabe Laplacian; this recovers an earlier result of K. Okikiolu.

\section{CR geometry}

Now we turn to the CR case.  Many of the results outlined in Section 2 have natural analogs in CR/complex setting.  However, we here put weight on the parts that are specific to the CR case  and will omit many of the fundamental results, for which we refer to
\cite{HPT} and \cite{GG}.  To simplify the exposition, we only consider the case of 
strictly pseudoconvex domains $\bC^{n+1}$ and the CR structure on the boundary.  More general formulation is given in \cite{H4}, in which Lee's pseudo-Einstein condition plays essential role.

\subsection{The ambient metric for strictly pseudoconvex domains}
Let $\Omega\subset\bC^{n+1}$ be a bounded strictly pseudoconvex domain with smooth boundary
$M=\pa\Omega$.
The CR structure on $M$ is given by 
$T^{1,0}=T^{1,0}\bC^{n+1}\cap \bC T M$, a rank $n$
complex subbundle of $\bC TM$.   It is {\em integrable} in the sense that 
$$
[\Gamma(T^{1,0}),\Gamma(T^{1,0})]\subset \Gamma(T^{1,0}).
$$ 
Take a $C^\infty$ defining function $\rho$ of $\Omega$ which is positive in $\Omega$, i.e., $\rho\in C^\infty(\bC^{n+1})$, $\Omega=\{\rho>0\}$ and $d\rho\ne0$ on $M$. Then we can define three hermitian metrics:
\begin{itemize}

\item {\bf The Levi metric}: for each $p\in M$,  a hermitian form  on $T_p^{1,0}$ is given
by
$$
L_\rho(Z,W)=-\pa\conj\pa\rho (Z,\conj W),
\qquad Z,W\in T_p^{1,0},
$$
which is positive  by the definition the strictly pseudoconvexity of $\Omega$.

\item {\bf Complete K\"ahler metric}: on the domain $\Omega$,  the real $(1,1)$-form
$$
g_+[\rho]=-\sqrt{-1}\pa\conj\pa\log\rho
$$
gives a complete K\"ahler metric near the boundary.

\item {\bf Lorentz-K\"ahler metric}: let $\rho_\sharp\colon\bC^*\times\bC^{n+1}\to \bR$ be the function 
$
\rho_\sharp(z_0,z)=|z_0|^2\rho(z)$, $(z_0,z)\in\bC^*\times\bC^{n+1}.
$
Then 
$$
\wt g[\rho]=-\sqrt{-1}\pa\conj\pa\rho_\sharp
$$
gives a Lorentz-K\"ahler metric on $\bC^*\times\bC^{n+1}$
near the hypersurface $\rho_\sharp=0$.
\end{itemize}
\medskip
In the case $\rho(z)=1-|z|^2$, $\Omega$ is the unit ball $B^{2n+2}\subset\bC^{n+1}$,
$$
(g_+[\rho])_{i\conj j}=
\rho(z)^{-1}\delta_{i\conj j}
-\rho(z)^{-2}\conj z_i z_j 
$$
is the complex hyperbolic metric (or the Bergman metric) and
$$
\wt g[\rho]=\sqrt{-1}\big(-d\zeta_0\wedge d\conj\zeta_0+
d\zeta_1\wedge d\conj\zeta_1+\cdots
+d\zeta_{n+1}\wedge d\conj\zeta_{n+1}\big)
$$
is a flat Lorentz-K\"ahler metric, where 
$\zeta$ is the coordinate system given by
$$
\zeta_0=z_0,\quad \zeta_j=z_0z_j,\quad j=1,\dots,n+1.
$$
As in the conformal case, the special unitary group $G=SU(n+1,1)$ for this metric acts on 
$\bC^{n+2}$ as linear transformations in $\zeta$ and induces isometries on $(B^{2n+2},g_+)$
and CR diffeomorphisms on $S^{2n+1}=\pa B^{2n+2}$, i.e., the diffeomorphisms that preserve the subbundle $T^{1,0}$.

For a general strictly pseudoconvex domain, we fix the defining function $\rho$
by imposing a complex Monge-Amp\`ere equation:
$$
\calJ_z[\rho]=1\qquad \text{on} \quad \Omega,
$$
where 
$$
\calJ_z[\rho]=(-1)^{n}\det\begin{pmatrix}\rho & \pa_j \rho\\ \pa_{\conj k}\rho & \pa_{j\conj k}\rho
\end{pmatrix}_{j,k=1,\dots,n+1}.
$$
The unique existence of the solution has been proved by S. Y. Cheng and S. T. Yau.  For such $\rho$, one has
$$
\begin{aligned}
\Ric[g_+]&=-(n+1)g_+\quad& \text{ on } \quad\rho>0,\\
\Ric[\,\wt g\,]&=0\quad&\text{on}\quad \rho_\sharp>0.
\end{aligned}
$$
However,  in general, the exact solution has weak singularity at the boundary.  We thus use
 the best approximate solution  $r\in C^{\infty}(\bC^{n+1})$ constructed by Fefferman \cite{FMA}.
Recall  that there is a smooth defining function $r$ such that
\begin{equation}\label{JeqC}
\calJ_z[r]=1+\eta\,r^{n+2}
\end{equation}
for an $\eta\in C^{\infty}(\bC^{n+1})$.
Such an $r$ is unique modulo $O(r^{n+3})$ and is called {\em Fefferman's defining function} of $\Omega$.
It is also important to note that
$$
\calO=\eta|_{M}\in C^{\infty}(M)
$$
is independent of the choice $r$ and
is called the {\em obstruction function}, as $\calO=0$  if and only if the Cheng-Yau solution is smooth up to the boundary.
We define the {\em ambient metric} to be $\wt g[r]$, which is not Ricci-flat but $\Ric[\,\wt g\,]=O(r^{n})$
holds.

The operator $\calJ_{z}$ depends on the choice of coordinates $z$ and so does $r$.
However, there is a simple transformation rule under the coordinate changes.  If $\wh z=\Phi(z)$ is another holomorphic coordinate system,
$$
\wh r=e^{-2\Re \varphi(z)} r 
$$
gives Fefferman's solution in $\wh z$, where 
$\varphi(z)=(\det\Phi'(z))^{1/(n+2)}$ is the power of holomorphic Jacobian.
Thus the map 
$$
\Phi_{\sharp}\colon (z_{0},z)\mapsto (\wh z_{0},\wh z)=(z_{0}\varphi(z), \Phi(z))
$$
 gives an isometry between $\wt g[r]$ and $\wt g[\,\wh r\,]$. In other word, we can say that the ambient metric is naturally defined on an $(n+2)$-nd root of the canonical bundle.

In the following, we regard the domain $\Omega\subset\bC^{n+1}$ as a complex manifold and consider the family of Fefferman's defining functions, each of them corresponds to a  choice of  coordinates.
For any defining functions $r$ and $\wh r$ in the family, we have 
$$
\wh r=e^{\Upsilon}r\quad\overline{\Omega}\mod O(r^{n+3})
$$
for a pluriharmonic function $\Upsilon$ on $\overline\Omega$.
If we define a contact form by 
$$
\th[r]=\frac{i}2(\pa-\conj\pa)r\big|_{TM},
$$
 then we have
$$
\wh\th=e^{\Upsilon}\th\quad\text{on}\quad{M}.
$$
Therefore the family of defining functions 
(or corresponding contact forms) can be seen as an analogy of conformal structure.
Important fact here is that the scaling is parametrized not by $C^{\infty}(M)$ but by the boundary values of pluriharmonic functions, which are called {\em CR pluriharmonic functions}.

\subsection{$Q$-prime curvature}
We  use the ambient metric $\wt g=\wt g[r]$  to construct CR invariant differential operators. Let $\wt\Delta$ be the Laplacian of $\wt g$. Then, for an integer $2m\in[-n,0]$ and a function $f\in C^{\infty}(\overline \Omega)$, 
$$
\Big(\wt\Delta^{n+2m+1} |z_0|^{2m} f\Big)\Big|_{\{1\}\times M}\in C^{\infty}(M)
$$ 
is shown to depend only on the boundary value of $f$ and gives a differential operator
$$
P_{n+2m+1}\colon C^{\infty}(M)\to C^{\infty}(M).
$$
While this definition depends on the choice of  $r$, we can say that this is
{\em CR invariant } in the sense that if $\wh r=e^{\Upsilon} r$ for a pluriharmonic function $\Upsilon$, then
$$
\wh P_{n+2m+1}(e^{m\Upsilon}f) =e^{(-n-m-1)\Upsilon} P_{n+2m+1}f.
$$
We will use the density notation and write this transformation law as
$$
P_{n+2m+1}\colon\calE(m)\to\calE(-n-m-1).
$$
The case $m=0$ has special importance as we have  $\calE(-n-1)=\Gamma(\wedge^{2n+1}T^{*}M)$ so that
$$
P_{n+1}\colon C^{\infty}(M)\to\Gamma(\wedge^{2n+1}T^{*}M).
$$

Now we recall Branson's idea of defining $Q$-curvature from these invariant operators.
Consider the $0$th order term of
$P_{n+1}$ for higher dimensions $\bC^{N+1}$ and take the ``limit as $N\to n$'' after factoring out $(N-n)$. This gives a formal definition of $Q$-curvature:
$$
Q_{n+1}=\lim_{N\to n}\frac{1}{N-n}\Big(\wt\Delta^{n+1} |z_0|^{2(n-N)}\Big)\Big|_{\{1\}\times M^{2N+1}}.
$$
We can justify this limit by considering Taylor expansion in $N-n$:
$$
|z_0|^{2(n-N)}=\sum_{k=0}^\infty\frac{(N-n)^k}{k!}
(-\log |z_0|^2)^k.
$$
Applying $\wt\Delta^{n+1}$ on $\bC^*\times\bC^{n+1}$ to  the both sides gives
$$
\big(\wt\Delta^{n+1}|z_{0}|^{2(n-N)}\big)|_{\{1\}\times M}=\sum_{k=0}^\infty\frac{(N-n)^k}{k}
Q^{(k)},
$$ 
where
$$
Q^{(k)}=\wt\Delta^{n+1}(-\log |z_{0}|^2)^k\big|_{\{1\}\times M}.
$$
While the expansion does not have clear meaning, the coefficients $Q^{(k)}$ are standard quantities defined on $M$ of dimension $2n+1$.
Clearly, $Q^{(0)}=0$.
In the conformal case, the second term $Q^{(1)}$ gives the $Q$-curvature, where
$\log|z_{0}|^{2}$ is replaced by $\log t^{2}$.  However, $Q^{(1)}=0$ because $\log |z_0|^2$ is pluriharmonic and $\wt g$ is K\"ahler.
Hence the leading term of the expansion is $Q^{(2)}$,
which we define to be the {\em Q-prime curvature} and denote by $Q'$.

The definition of $Q'$ depends on the choice of $r$ and is not a CR invariant.  
If $\wh r=e^\up r$, where $\up$ is pluriharmonic, then
$$
\wh Q'\ =Q'+2P'\up+P_{n+1}(\up^2).
$$
Here $P'$ is a differential operator defined on the space of CR pluriharmonic functions 
by
$$
P'f=-\wt\Delta^{n+1}(\wf\log |z_0|^2)\big|_{\{1\}\times M},
$$
which we call the {\em $P$-prime operator}.  Here $\wt f$ denotes the pluriharmonic extension of $f$.
Again, $P'$ is not a CR invariant operator but satisfies the transformation law:
$$
\wh P'f=P'f+P_{n+1}(\up f).
$$
A crucial fact is that $P'$ and $P_{n+1}$ are formally self-adjoint and $P'1=P_{n+1}1=0$.  It follows that the  {\em total $Q$-prime curvature}
$$
\conj{Q}'(M)=\int_M Q'\th\wedge(d\th)^{n}
$$
is a CR invariant of $M$, i.e., it is independent of the choice of $r$.

\subsection{Explicit formulas in dimensions $3$ and $5$}
In the case $M$ has dimension 3, we can explicitly write down $P'$ and $Q'$ in terms of Tanaka-Webster connection $\nabla$ (analogous to the Levi-Civita connection, for each choice of a contact form or the Levi metric, one can define a canonical connection of $TM$). With respect to the contact form $\th=\th[r]$ for $r$ given as above, we have 
\begin{align}\label{PP3}
P'f&=\Delta_b^2f-\Re\nabla^1(\Scal \nabla_1f-2\sqrt{-1}A_{11}\nabla^1f),
\\
\label{Q4form}
Q'&=\frac12\Delta_b \Scal+\frac14 \Scal^2-|A|^2.
\end{align}
Here $\Delta_b$ is the sub-Laplacian, and $\Scal$, $A_{11}$ are respectively the scalar curvature and torsion of the connection; $|A|^{2}=A_{11}A^{11}$ is the squared norm of the torsion.  We are still using the Einstein convention but, since $T^{1,0}$ has rank one, we only have index $1$.

These formulas were first given by J. Case and P. Yang  prior to the general definition in the previous subsection.
Their aim was to give a CR analogue of Gursky's sphere theorem in 4-dimensional conformal geometry.
While we cannot go into the details, let us recall their main theorem.

\begin{theorem} [Case-Yang \cite{CY}] 
Let $(M,H,J)$ be a compact 3-dimensional CR manifold with a pseudo-Einstein contact form.
Assume that $P_{3}$ is nonnegative and that CR Yamabe constant is nonnegative.
Then
$$
\overline Q'(J)\le \overline Q'(J_{0})
$$
and the equality holds if and only if $(M,H,J)$ is CR equivalent to the standard sphere $(S^{3},H_{0},J_{0})$.
\end{theorem}

This is a deep result of geometric analysis; the proof is based on the CR positive mass theorem of 
J.-H. Cheng, A. Malchiodi, and P. Yang.

From the explicit formula of $Q'$ in 3-dimensions, we can see that $\overline Q'$ agrees with the Burns-Epstein invariant $\mu(M)$, \cite{BE1}, up to a universal constant:
$$
\conj{Q}'(M)=-4\pi^2\mu(M).
$$
From this fact we can also obtain the renormalized Gauss-Bennet formula for $\Omega\subset\bC^2$:
\begin{equation}\label{GBthm}
\int_\Omega c_2(B)=\chi(\Omega)-\frac{1}{4\pi^2}\conj{Q}'(M),
\end{equation}
where $c_{2}$ is the second Chern form for the Bochner tensor $B$ of $g_{+}$, the trace-free part of the K\"ahler curvature tensor of $g_{+}$.

For higher dimensions, such equality does not hold in general.
To state it precisely, let us recall a result of T. Marugame \cite{Maru}, which 
improved the renormalized Gauss-Bonnet formula of Burns-Epstein \cite{BE2}. 
For $\Omega\subset\bC^{n+1}$, he found a transgression formula that gives an invariant polynomial $\Pi(R,A)$ in the curvature $R$ and torsion $A$ of the Webster-Tanaka connection for $\th$ such that
$$
\int_\Omega c_{n+1}(B)=\chi(\Omega)-\int_{M}\Pi\cdot\th\wedge(d\th)^{n},
$$
where $c_{n+1}(B)$ is the $(n+1)$-st Chern form for the Bochner tensor of $g_{+}$ on $\Omega$.
When $n=2$, we have
\[
-(4\pi)^3\Pi= 
\frac{1}{27}\Scal^3-4R_{a\ol{c}b\ol{d}}A^{ab}A^{\ol{c}\ol{d}}+\frac{1}{3}|S|^2\Scal.
\]
Here  $|S|^2$ is the squared norm of the Chern-Moser tensor $S_{a\ol{c}b\ol{d}}$,
which is the trace-free part of the Tanaka--Webster curvature $R_{a\ol{c}b\ol{d}}$. 
 The lower indices $a,b$ (resp.~$\conj c,\conj d$) run through $1,2,\dots,n$  (resp.~$\conj 1,\conj 2,\dots,\conj n$) and correspond to $(T^{1,0})^{*}$  (resp.~$(T^{1,0})^{*}$).  
  Analogous to the Weyl curvature in the conformal case, $S_{a\ol{b}c\ol{d}}=0$ if and only if $M$ is {\em spherical}, i.e.,
 locally CR equivalent the sphere.
 With this $\Pi$, we can write $\overline Q'$ as
\[
\ol{Q}^{\prime}=-\int_M \Bigl(
(4\pi)^3\Pi+\frac{1}{3}|S|^2\Scal +4|\nabla A|^{2}\Bigr)\th\wedge(d\th)^2.
\]
Thus we obtain the following

\begin{theorem}[
 \cite{HMM}] Let $\Omega\subset \bC^{3}$ be a strictly pseudoconvex domain.
If the scalar curvature of the Webster-Tanaka connection for $\th[r]$ is positive almost everywhere, then
$$
(4\pi)^3\overline Q'\le \chi(\Omega)-\int_\Omega c_{3}(B).
$$
The equality holds only if $M=\pa\Omega$ is spherical.
\end{theorem}

Note that the assumption on $\Omega$ holds for domains that are sufficiently close to the ball in $\bC^3$. So there are many examples for which 
$ -(4\pi)^3\overline Q'$ and  $\int_{M}\Pi$ are different.

\subsection{Volume renormalization and variational formula}
We next consider the volume renormalization of strictly pseudoconvex domains.
While $M$ has odd dimensions, 
CR geometry is analogous to even dimensional conformal geometry.
Hence the situation is a little bit  complicated.
As before, let  $g_{+}=-\sqrt{-1}\pa\conj\pa\log r$ be the complete K\"ahler metric on $\Omega$ defined from Fefferman's defining function $r$. 

\begin{theorem}[\cite{H4,HMM}]\label{renormalized-volume-theorem}  Let $\Omega\subset\bC^{n+1}$ be a
strictly pseudoconvex domain and $r$ be Fefferman's defining function. 
 Then the integrals over the subdomains $\{r>\epsilon\}$ admit expansions as $\epsilon\to +0$:
\begin{align}\label{renormalization1}
\int_{r>\epsilon}\abs{d\log r}_{g_+}^2d\vol_{g_+}
	&=a_0\epsilon^{-n-1}+\dots+a_n\epsilon^{-1}+k_n\overline Q'\log{\epsilon}+O(1),
\\
\int_{r>\epsilon}d\vol_{g_+}&=b_0\epsilon^{-n-1}+\dots+b_n\epsilon^{-1}+k_n'\overline Q'+o(1),
\end{align}
where $a_j,b_j$ are constants given by integrals over $M=\pa\Omega$ of some local invariants of the CR structure of $M$ and $\th=\th[r]$, and $k_n,k_n'$ are non-zero universal constants depending only on the dimension.
\end{theorem}

The first formula is an analogy of the even dimensional conformal case,  while the second formula says that the renormalized volume is a CR invariant, which corresponds to the odd dimensional conformal case.

Using \eqref{renormalization1}, we can compute the variation of $\overline Q'$ under the perturbation of domains.
Let $\{\Omega_t\}_{t\in\bR}$ be a smooth family of strictly pseudoconvex domains in $\bC^{n+1}$ in the sense that there is a $C^\infty$ function $\rho_t(z)$ of $(t,z)\in\bR\times\bC^{n+1}$ such that $\Omega_t=\{z\in\bC^{n+1}:\rho_t(z)>0\}$ and $d_z\rho_t\ne0$ on $\pa\Omega_t$.
Solving the Monge-Amp\`ere equation for each $t$, one may assume that $\rho_{t}$ is Fefferman's defining function for each fixed $t$.
On the boundary $M_{t}=\pa\Omega_{t}$, $\th_t=\th[\rho_t]$ gives a natural contact form. 

\begin{theorem}[\cite{HMM}]\label{first_variation_theorem}
Let $\{\Omega_t\}_{t\in\bR}$ be a smooth family of strictly pseudoconvex domains in $\bC^{n+1}$.  Then the total $Q'$-curvature $\overline Q'(M_{t})$ of $M_t$ satisfies
\begin{equation}\label{QPvar}
\frac{d}{d t}\Big|_{t=0}\overline Q'(M_{t})=2 \int_{M_0} \dot\rho\,\calO\,\th_{0}\wedge(d\th_{0})^{n},
\end{equation}
where 
$
\dot\rho(z)={d\rho_{t}}/{d t}|_{t=0}
$
and $\calO$ is the obstruction function of $\rho_0$.
\end{theorem}

\subsection{$Q$-curvature for partially integrable CR structures}
As we have seen, the $Q$-curvature vanishes for the boundary of a domain in $\bC^{n+1}$. However, if we consider abstract CR structure which may not be embeddable, the total $Q$-curvature becomes non-trivial and  has natural variational formula.
We here recall a result of Y. Matsumoto \cite{Ma2}.

 Let $\th$ be a contact form on a manifold $M$ of dimension $2n+1$, that is, $\th$ is a real one form satisfying $\th\wedge(d\th)^{n}\ne0$.  
An abstract CR structure is a complex structure $J$ on the contact distribution $H=\ker\th\subset TM$; we denote the  $\pm\sqrt{-1}$-eigenspace decomposition by
 $\bC H=T^{1,0}\oplus T^{0,1}$.  We assume that $J$ is {\em partially integrable} in the sense that
$$
[\Gamma(T^{1,0}),\Gamma(T^{1,0})]\subset \Gamma(T^{1,0}\oplus T^{0,1}).
$$
This enables us to define the Levi metric $L_{\th}$ on $H$ by $d\th(X,JY)$;  we  assume that it is positive definite. 
A choice of contact form $\th$ gives a decomposition of cotangent bundle
$\bC T^{*}M=\bC \th\oplus (T^{1,0})^{*}\oplus (T^{0,1})^{*}$ such that the corresponding coframe
$\th,\th^{a},\th^{\conj a}$ gives $d\th=\sqrt{-1}h_{a\conj b}\th^{a}\wedge\th^{\conj b}$.

Let $(M,H,J)$ be a partially integrable CR manifold with a contact form $\th$.
An {\em asymptotically complex hyperbolic (ACH) metric} is a Riemannian metric on $M\times (0,1)$ with the following asymptotic expansion
$$
g_{+}=\frac{1}{4\rho^2}d\rho^2+\frac1{\rho^2}\th^{2}+O(\rho^{-1}),\qquad
g_{+}|_{H}=\frac{L_{\th}}{\rho}+O(1).
$$
Here $\th$ and $L_{\th}$ are identified with their pullbacks by the projection
$M\times(0,1)\to M$.  We also assume that $g_+$ is {\em smooth}  in the sense that $\rho^2g_+$ is $C^\infty$ on $M\times[0,1)$.

As in the case of Poincare-Einstein metric, we consider the best approximate solution to the Einstein equation.  It is shown that there is a smooth ACH metric that satisfies
\begin{equation}\label{ACHEq}
\Ric(r_+)+\frac12(n+2)g_{+}=\rho^n E
\end{equation}
for a symmetric 2-tensor  $E$ which is $C^\infty$ on $M\times [0,1)$
and 
$$
E |_{\rho=0}=2\Re(\calO_{ab}\th^a\otimes\theta^b)
\mod\th, d\rho.
$$
The {\em obstruction tensor} in this setting is defined to be
$$
\calO_{ab}\in \Gamma(S^{2} (T^{1,0}M)^*).
$$
It is CR invariant, i.e., 
 $\wh\calO_{\alpha\beta}=e^{-n\Upsilon}\calO_{\alpha\beta}$ holds
 under the scaling $\wh\th=e^{\Upsilon}\th$.
If $\calO_{\alpha\beta}=0$, we can find a smooth ACH  metric that satisfies Einstein equation modulo $O(\rho^\infty)$.
It is important to note that  $\calO_{\alpha\beta}=0$ if $T^{1,0}$ is integrable. 
  
Fixing a smooth ACH metric $g_{+}$ satisfying \eqref{ACHEq},
 we now define $Q$-curvature for $\th$ on $(M,H, J)$.
Let $\Delta_{+}$ be the Laplacian of $g_+$. Then there are functions $A,B\in C^{\infty}(M\times[0,1) )$ such that
$$
\Delta_{+}(\log\rho+A+B\rho^{n+1}\log\rho)=n+1+O(\rho^{\infty}).
$$
The $Q$-curvature is now defined by
$$
B|_{\rho=0}=\frac{(-1)^{n}}{n!(n+1)!}\,Q.
$$
One can show that $Q$ satisfies the required transformation law
$$
\wh Q=Q+P_{n+1}\Upsilon,
\qquad\wh\th=e^{\Upsilon}\th, 
$$
where $\Upsilon\in C^{\infty}(M)$ and $P_{n+1}$ is a self-adjoint CR invariant differential operator of order $2n+2$ without constant term.
It follows that 
$$
\overline Q=\int_M Q\th\wedge(d\th)^n
$$
is a CR invariant. 
For integrable CR structures, this definition of $Q$ agrees with the one given in \cite{FH}
via the ambient metric
(in this case, we can also say that CR $Q$ is the pushforward of the conformal $Q$ of the Fefferman space $S^{1}\times M$).

To state the variational formula of $\overline Q$, we recall the deformation of (partially integrable) CR structures.
Take a frame $Z_{a}$ of $T^{1,0}$ and set   $Z_{\conj a}=\conj{Z_a}\in T^{0,1}$.
Then we may define another CR structure $\wh J$ by the frame of $\wh T^{0,1}$,
$$
\wh Z_{\conj a}=Z_{\conj a}+\varphi_{\conj a}{}^{b}Z_{ b}.
$$
It is partially integrable if $\varphi_{\conj a\conj b}$ is symmetric, where the index $b$ is lowered by using the Levi metric $h_{a\conj b}$.  Thus partially integrable CR structure nearby $T^{1,0}$ is parametrized by a symmetric two tensor $\varphi_{\conj a\conj b}\in\Gamma(S^{2}(T^{0,1})^{*})$.

\begin{theorem}[Matsumoto \cite{Ma2}]
Let  $\{J_{t}\}_{t\in\bR}$ be a one parameter family of partially integrable CR structures parameterized by $\varphi_{\conj a\conj b}^{(t)}$.  Then the total $Q$-curvature  $\overline Q(J_{t})$ of $(M,H,J_{t})$ satisfies
$$
\frac{d}{dt}\Big|_{t=0}\overline Q(J_{t})=
(-1)^{n}c_n\int \Re \calO_{ab}\dot\varphi^{ab}\th\wedge(d\th)^{n}
$$
with a universal constant $c_n>0$.   Here $\calO_{ab}$ is the obstruction tensor for
$(M, H,J_{0})$ in the scale $\th$ and
$
\dot\varphi_{\conj a\conj b}={d}\varphi_{\conj a\conj b}^{(t)}/{dt}|_{t=0}.
$
\end{theorem}

In particular, if all $J_{t}$ are integrable, we have $\calO_{ab}=0$ and
$\overline Q(J_{t})$ is constant.  

\subsection{Deformation complex of CR structures}

We have obtained the variational formulas of $\overline Q$ and $\overline Q'$.  To derive geometric consequences from them, we shall recall the deformation complex of CR structures.

We will use the Young diagram to denote the symmetries of tensor bundle.  The symmetric product of $(T^{0,1})^*$ is now denoted by $\YDsym(T^{0,1})^*
$.
The integrability of $J$ is equivalent to the vanishing of the Nijenhuis tensor $N_{\conj a\conj b\conj c}$, which has the symmetry $\YDN$; the linearization of $N_{\conj a\conj b\conj c}$ gives the map
$$
D_1^-\colon 
\Gamma\Big(\,
\setlength{\unitlength}{0.25mm}
\YDsym(T^{0,1})^*\Big)
\to
\Gamma\Big(\,
\YDN(T^{0,1})^*
\Big),  \qquad D_1^- \varphi_{\conj a\conj b}=\nabla_{[\conj c}\varphi_{\conj b]\conj a}.
$$
Let $\calE=C^{\infty}(S^{2n+1},\bC)$ and define
$$
\begin{aligned}
D_0^-\colon&\calE\to
\Gamma\Big(\,
\YDsym (T^{0,1})^{*}
\Big),  \quad &&D_0^- f=\nabla_{\conj a\conj b} f,
\\
D_0^+\colon&\calE\to\Gamma\Big(\YDsym (T^{1,0})^{*}\Big),
 && D_0^+ f=\nabla_{ab} f.
\end{aligned}
$$
These maps give a complex
$$
\calE\xrightarrow{\ D_{0}^{-} \ }
\Gamma\Big(\,
\YDsym(T^{0,1})^{*}
\Big)\xrightarrow{\ D_{1}^{-} \ }
\Gamma\Big(\,
\YDN(T^{0,1})^{*}
\Big).
$$
This is the beginning of Kuranishi's deformation complex of CR structures in the form later improved by Akahori-Garfield-Lee \cite{AGL}.  The cohomology of this complex describes the moduli of the deformations of isolated singularities.  In this setting, CR manifolds that bound the same singularity are identified; this equivalence is given infinitesimally by the image of $D_{0}^{-}$.

To study the deformation of partially integrable CR structures, we need the full deformation complex which is given as  the BGG complex of the adjoint representation $\mathfrak{su}(n+1,1)$. The Kuranishi complex is contained in the BGG complex as an edge.  
On the 5-dimensional sphere $S^{5}\subset\bC^3$, the deformation complex is given by
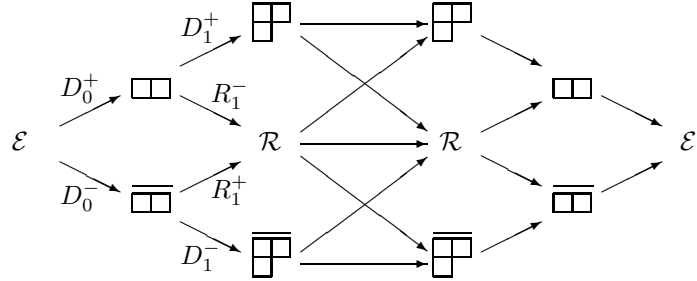
\begin{equation*}
\setlength{\unitlength}{0.8mm}
\begin{picture}(120,50)(0,-25)
\put(2,-1){$\calE$}
\put(10,8){$D^+_0$}
\put(10,-10){$D^-_0$}
\put(10,3){\vector(2,1){10}}
\put(10,-3){\vector(2,-1){10}}
\put(22,8){$\YDsym$}
\put(22,-11){$\overline\YDsym$}
\put(30,18){$D^+_1$}
\put(35,7){$R^-_1$}
\put(35,-9){$R^+_1$}
\put(30,-20){$D^-_1$}
\put(30,13){\vector(2,1){10}}
\put(30,8){\vector(2,-1){10}}
\put(30,-8){\vector(2,1){10}}
\put(30,-13){\vector(2,-1){10}}
\put(42,19){$\YDN$}
\put(43,-1){$\calR$}
\put(42,-22){$\overline\YDND$}
\put(50,20){\vector(2,0){21}}
\put(50,18){\vector(4,-3){21}}
\put(50,2){\vector(4,3){21}}
\put(50,0){\vector(2,0){21}}
\put(50,-2){\vector(4,-3){21}}
\put(50,-18){\vector(4,3){21}}
\put(50,-20){\vector(2,0){21}}
\put(72,19){$\YDN$}
\put(73,-1){$\calR$}
\put(72,-22){$\overline\YDND$}
\put(80,18){\vector(2,-1){10}}
\put(80,2){\vector(2,1){10}}
\put(80,-2){\vector(2,-1){10}}
\put(80,-18){\vector(2,1){10}}
\put(92,8){$\YDsym$}
\put(92,-11){$\overline\YDsym$}
\put(100,-8){\vector(2,1){10}}
\put(100,8){\vector(2,-1){10}}
\put(113,-1){$\calE$}
\end{picture}
\end{equation*}
Here we have simplified the notation by omitting $\Gamma$ and  $(T^{1,0})^{*}$ or $(T^{0,1})^{*}$; the overline means that $(T^{0,1})^{*}$ is omitted. 
$\calR$ denotes space of the sections of the trace-free tensors $\varphi_{ab\conj c\conj d}$
with symmetry $\YDsym\otimes\conj\YDsym$ and $R_{1}^{-}$ is given by
$$
R_{1}^{-}(\varphi_{ a b})=\text{trace-free part of }\nabla_{ \conj c\conj d}\varphi_{ a b}.
$$
Each arrow is a CR invariant differential operator if we properly put density weight on each bundle. Moreover,  it is known that there is exactly one CR invariant operator for each arrow.
(Recall that $G=SU(n+1,1)$ acts on $S^{2n+1}$ as CR automorphisms. Hence, in this setting, CR invariant operators are $G$-invariant operators and vice versa.)

For higher dimensions, we have a diagram of length $2n+2$.  It is similar to the type decomposition of the de Rham complex, but it has more maps in the middle.  The part we need here is the following:
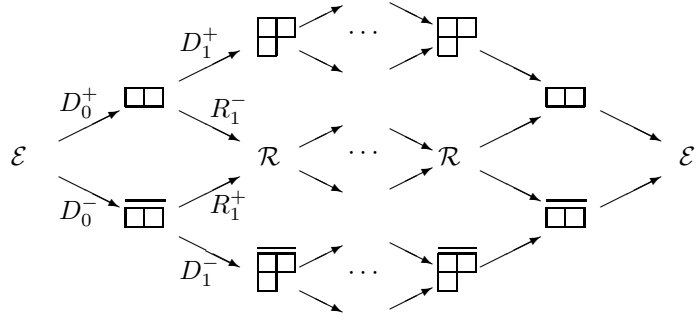
\begin{equation*}
\setlength{\unitlength}{0.8mm}
\begin{picture}(115,55)(0,-25)
\put(2,-1){$\calE$}
\put(10,8){$D^+_0$}
\put(10,-10){$D^-_0$}
\put(10,3){\vector(2,1){10}}
\put(10,-3){\vector(2,-1){10}}
\put(21,9){$\YDsym$}
\put(21,-11){$\overline\YDsym$}
\put(30,18){$D^+_1$}
\put(35,7){$R^-_1$}
\put(35,-9){$R^+_1$}
\put(30,-20){$D^-_1$}
\put(30,13){\vector(2,1){10}}
\put(30,8){\vector(2,-1){10}}
\put(30,-8){\vector(2,1){10}}
\put(30,-13){\vector(2,-1){10}}
\put(43,19){$\YDN$}
\put(43,-1){$\calR$}
\put(43,-22){$\overline\YDND$}
\put(50,22){\vector(2,1){7}}
\put(50,18){\vector(2,-1){7}}
\put(50,2){\vector(2,1){7}}
\put(50,-2){\vector(2,-1){7}}
\put(50,-18){\vector(2,1){7}}
\put(50,-22){\vector(2,-1){7}}
\put(58,20){$\cdots$}
\put(58,0){$\cdots$}
\put(58,-20){$\cdots$}
\put(65,26){\vector(2,-1){7}}
\put(65,15){\vector(2,1){7}}
\put(65,5){\vector(2,-1){7}}
\put(65,-5){\vector(2,1){7}}
\put(65,-14){\vector(2,-1){7}}
\put(65,-25){\vector(2,1){7}}
\put(73,19){$\YDN$}
\put(73,-1){$\calR$}
\put(73,-22){$\overline\YDND$}
\put(80,18){\vector(2,-1){10}}
\put(80,2){\vector(2,1){10}}
\put(80,-2){\vector(2,-1){10}}
\put(80,-18){\vector(2,1){10}}
\put(91,9){$\YDsym$}
\put(91,-11){$\overline\YDsym$}
\put(100,-8){\vector(2,1){10}}
\put(100,8){\vector(2,-1){10}}
\put(113,-1){$\calE$}
\end{picture}
\end{equation*}
If we write $f=u+\sqrt{-1}v$ for real valued $u$ and $v$, then $D_0^{-} \sqrt{-1}v$ is infinitesimally given by a pullback of the CR structure by a contact diffeomorphism, while $D_0^{-} u$ is the first variation of Kuranishi wiggle (perturbations of $S^{2n+1}$ within $\bC^{n+1}$).

We set $\calQ=\bC T/(T^{1,0}\oplus T^{0,1})$ and denote the space of sections of $\calQ$ by $\calE(1)$.  For each choice of a contact form
we have an identification $\calE(1)\cong \calE$. In the deformation complex,
the first $\calE$ is $\calE(1)$ and the last one is $\calE(-n-2)=\Gamma(\calQ^{-n-2})$.
The symmetric two tensor on the left side $\YDsym$ has weight $\YDsym\otimes\calQ$ and the one on the right side has weight $\YDsym\otimes\calQ^{-n}$;
we set  
$$
\calE_\YDsymsc(1)=\Gamma(\YDsym\otimes\calQ)
\quad\text{ and }\quad
\calE_\YDsymsc(-n)
=\Gamma(\YDsym\otimes\calQ^{-n}).
$$ 

Take a family of partially integrable CR structures
$J_{t}$ given by
$\varphi^{(t)}_{ab}$ (the conjugate of $\varphi^{(t)}_{\conj a\conj b}$) and set 
 $$
\dot \varphi_{ab}=\frac{d}{dt}\Big|_{t=0}\varphi_{ab}^{(t)}\in\calE_\YDsymsc(1).
$$
Then the first variation of $\calO_{ab}(J_{t})$ at $t=0$ gives a CR invariant operator
of order $2n+2$:
$$
L_{n+1}\colon\calE_\YDsymsc(1)\to\calE_\YDsymsc(-n).
$$
For  the obstruction function $\calO$, we take its variation under an integrable deformation $\varphi_{ab}^{(t)}$ such that
$$
\frac{d}{dt}\Big|_{t=0}\varphi_{ab}^{(t)}=D_0^{+}f\in\calE_\YDsymsc(1)
\quad \text{for}\quad f\in\calE(1).
$$
 Then the first variation of $\calO(J_t)$ gives
a CR invariant operator of order $2n+6$:
$$
L_{n+3}\colon\calE(1)\to\calE(-n-2).
$$
These are the operators sending the  bundles on the left to the ones on the right with the same Young diagrams:
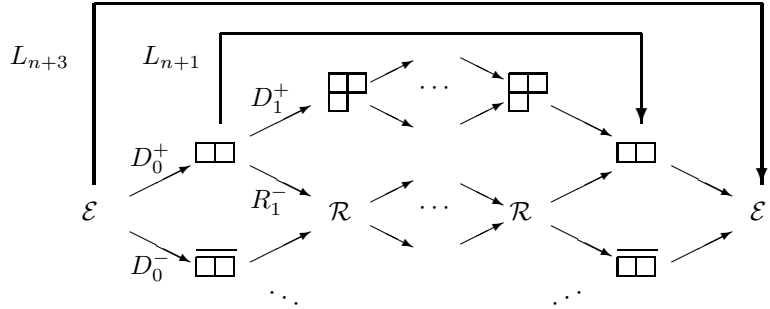
\begin{equation*}
\setlength{\unitlength}{0.8mm}
\begin{picture}(120,58)(0,-18)
{\thicklines
\put(-10,25){$L_{n+3}$}
\put(4,5){\line(0,1){30}}
\put(4,35){\line(1,0){111}}
\put(115,35){\vector(0,-1){30}}
\put(12,25){$L_{n+1}$}
\put(25,15){\line(0,1){15}}
\put(25,30){\line(1,0){70}}
\put(95,30){\vector(0,-1){15}}
}
\put(2,-1){$\calE$}
\put(10,8){$D^+_0$}
\put(10,-10){$D^-_0$}
\put(10,3){\vector(2,1){10}}
\put(10,-3){\vector(2,-1){10}}
\put(21,9){$\YDsym$}
\put(21,-10){$\overline\YDsym$}
\put(30,18){$D^+_1$}
\put(30,1){$R^-_1$}
\put(30,13){\vector(2,1){10}}
\put(30,8){\vector(2,-1){10}}
\put(30,-8){\vector(2,1){10}}
\put(33,-14){$\cdot$}
\put(35,-15){$\cdot$}
\put(37,-16){$\cdot$}
\put(43,19){$\YDN$}
\put(43,-1){$\calR$}
\put(50,22){\vector(2,1){7}}
\put(50,18){\vector(2,-1){7}}
\put(50,2){\vector(2,1){7}}
\put(50,-2){\vector(2,-1){7}}
\put(58,20){$\cdots$}
\put(58,0){$\cdots$}%
\put(65,26){\vector(2,-1){7}}
\put(65,15){\vector(2,1){7}}
\put(65,5){\vector(2,-1){7}}
\put(65,-5){\vector(2,1){7}}
\put(73,19){$\YDN$}
\put(73,-1){$\calR$}
\put(80,18){\vector(2,-1){10}}
\put(80,2){\vector(2,1){10}}
\put(80,-2){\vector(2,-1){10}}
\put(84,-14){$\cdot$}
\put(82,-15){$\cdot$}
\put(80,-16){$\cdot$}
\put(91,9){$\YDsym$}
\put(91,-10){$\overline\YDsym$}
\put(100,-8){\vector(2,1){10}}
\put(100,8){\vector(2,-1){10}}
\put(113,-1){$\calE$}
\end{picture}
\end{equation*}
Moreover, $L_{n+3}$ and $L_{n+1}$ are intertwiners between $G$-modules
$\calE(1)\to\calE(-n-2)$ and $\calE_\YDsymsc(1)\to\calE_\YDsymsc(-n)$; such maps are unique up to a constant multiple and the eigenvalues of these maps  can be explicitly computed by using representation theory.

\begin{theorem}[\cite{HMO}]\label{mainthem1}

\noindent
$(1)$ On the standard sphere $S^{2n+1}$
the operator
$L_{n+3}\colon \calE(1)\to\calE(-n-2)$ is semi-negative and
$$
\ker L_{n+3}=\ker D_0^++\ker D_0^-.
$$ 

\noindent
$(2)$ On the standard sphere $S^{2n+1}$, $n\ge2$, 
the operator
$L_{n+1}\colon \calE_\YDsymsc(1)\to\calE_\YDsymsc(-n)$ is semi-negative and 
$$
\ker L_{n+1}=\ker D_1^+ +\ker R_1^-.
$$

\end{theorem}

As a direct consequence, we have
\begin{theorem} [\cite{HMO}]\label{mainthm}
 
\noindent
$(1)$ Let $\Omega_{t}=\{\rho_{t}>0\}$ be a smooth family of strictly pseudoconvex domains in $\bC^{n+1}$ such that $\Omega_{0}$ is the unit ball. If  $\dot\rho=d\rho_{t}/{dt}|_{t=0}$
satisfies $\dot\rho|_{S^{2n+1}}\not\in\ker D_{0}^{+}+\ker D_{0}^{-}$, then
$$
\frac{d^{2}}{dt^{2}}\Big|_{t=0}\overline{Q}'(\pa\Omega_{t})<0.
$$
 
\noindent
$(2)$  Let $\{J_{t}\}_{t\in\bR}$ be a family of partially integrable CR structures such that
 $J_{0}$ is the standard one on $S^{2n+1}$, and $\dot\varphi_{ab}$ be the first variation of $J_{t}$ at $t=0$.  If $ \dot\varphi_{a b}\not\in\ker D_{1}^{+}+\ker R_{1}^{-}$, then
 $$
\frac{d^{2}}{dt^{2}}\Big|_{t=0}\overline{Q}(J_{t})<0.
$$
\end{theorem}

We shall explain the geometric meaning of the condition on the direction of the deformations.
If $n\ge2$,  we have 
$$
\ker D_0^++\ker D_0^-=\ker R_1^- D_0^+
$$
 and the composition
$R_{1}^{-}D_0^{+}$ is the linearization of the Chern-Moser tensor $S_{ab\conj c\conj d}$.
Thus  $\dot\rho|_{S^{2n+1}}\in\ker D_{0}^{+}+\ker D_{0}^{-}$ means that the family is spherical to the first order. The inequality states that if the family deformation is not trivial, then
$\overline Q'$ takes local maximum value at the standard sphere along the family.  

The similar argument can be applied to the case $n=1$; then $R_1^-$ is a 4-th order operator and
Chern-Moser tensor has 6 indices $S_{1111\conj 1\conj 1}$.  It is also not difficult to derive   Theorem \ref{mainthem1} (1), $n=1$, by a direct computation as was done in Burns-Epstein \cite{BE1}. J. H. Cheng and J. Lee further showed the following stronger result:

\begin{theorem} [Cheng-Lee \cite{CL}] For the CR structures on $S^3$ near the standard one,
 the Burns-Epstein invariant takes minimal value only for the standard sphere.
\end{theorem}

Since $\mu=-(2\pi)^{2}\overline{Q}'$, this theorem is consistent with the theorem above.
To prove the local minimality from the semi-positivity of the Hessian, they developed a slice theorem of the  moduli space of CR structures on $S^{3}$.

The geometric meaning of  Theorem \ref{mainthem1} (2) is still not clear; but it gives a insight to the partially integrable CR structures.
The subspace   $\ker D_{1}^{+}$ is the direction of integrable CR structures
and $\ker R_{2}^{-}$ is the direction with $S_{ab\conj c\conj d}=0$ (in the partially integrable case, $S_{ab\conj c\conj d}=0$ may not mean that the surface is spherical).
Thus the vanishing of  total $Q$ does not characterize integrable CR structures --- contrary to our initial hypothesis. However, this theorem suggests the existence of a natural class of partially integrable CR structures for which $\overline{Q}$ vanishes identically.

\end{document}